\documentclass[10pt,a4paper,oneside]{amsart}
\newcounter{commentcounter}

\usepackage{amsmath} % Lots of maths functionality
\usepackage{amssymb} % Maths symbols
\usepackage{amsthm} % Maths environments: \begin{proof}, etc.
\usepackage{stmaryrd} % [[ brackets
\usepackage[english]{babel} % Language and hyphenation.
\usepackage[font=small,justification=centering]{caption} % More flexibility for captioning figures
\usepackage[nodayofweek]{datetime}
\usepackage[shortlabels]{enumitem} % Change enumeration labelling with \begin{enumerate}[a)] etc.
\usepackage[T1]{fontenc} % For font encoding, to allow accents, copy & paste, inequality signs, etc. to all work nicely.
\usepackage[utf8]{inputenc} % To be loaded after fontenc, also for encoding.
\usepackage{ifthen} % For Dani's \begin{com} environment and \numberedtheorem
\usepackage{mathabx} % Contains \pnest symbol and more. Clashes with accents for \ring
\usepackage{mathtools} % Uses amsmath, fixes quirks and adds functionality.
\usepackage[dvipsnames]{xcolor} % Allow links to have colour. Needs to be before hyperref and tikz-cd
\usepackage[pdftex,  colorlinks=true]{hyperref} % Makes references and citations into links
    \hypersetup{urlcolor=RoyalBlue, linkcolor=RoyalBlue,  citecolor=black}
\usepackage{tikz-cd} % Commutative diagrams
\usepackage{xfrac} % Nicer quotients 
\usepackage[capitalize]{cleveref} % use \Cref{} for instead of X~\ref{} 
\makeatletter
\@namedef{subjclassname@1991}{Mathematical subject classification 1991}
\@namedef{subjclassname@2000}{Mathematical subject classification 2000}
\@namedef{subjclassname@2010}{Mathematical subject classification 2010}
\@namedef{subjclassname@2020}{Mathematical subject classification 2020}
\makeatother

\newtheorem{thm}{Theorem}[section]
\newtheorem{lemma}[thm]{Lemma}

\newtheorem{prop}[thm]{Proposition}

\newtheorem{thmx}{Theorem}
\newtheorem{corx}[thmx]{Corollary}

\theoremstyle{definition}
\newtheorem{defn}[thm]{Definition}
\newtheorem{remark}[thm]{Remark}

\theoremstyle{plain}

    \newtheoremstyle{TheoremNum}
        {\topsep}{\topsep} %%% space between body and thm
        {\itshape} %%% Thm body font
        {-0.25cm} %%% Indent amount (empty = no indent)
        {\bfseries} %%% Thm head font
        {.} %%% Punctuation after thm head
        { }  %%% Space after thm head
        {\thmname{#1}\thmnote{ \bfseries #3}}%%% Thm head spec
    \theoremstyle{TheoremNum}
    \newtheorem{duplicate}{}

\newcommand*{\claimproofname}{My proof}

\def\Z{\mathbb{Z}}

\usepackage{tikz}
\usetikzlibrary{arrows,quotes}
\tikzstyle{blackNode}=[fill=black, draw=black, shape=circle]

\title[Exotic subgroups of hyperbolic groups]{An infinite family of exotic subgroups of hyperbolic groups}
\author{Monika Kudlinska}
\address{Mathematical Institute, Andrew Wiles Building, Observatory Quarter, University of Oxford, Oxford OX2 6GG, UK}
\email{monika.kudlinska@maths.ox.ac.uk}
\date{\today}


\begin{document}

\maketitle

\begin{abstract}
    We construct the first known infinite family of quasi-isometry classes of subgroups of hyperbolic groups which are not hyperbolic and are of type $\mathrm{FP}(\mathbb{Q})$. We give a simple criterion for producing many non-hyperbolic subgroups of hyperbolic groups with strong finiteness properties. We also observe that local hyperbolicity and algebraic fibring are mutually exclusive in higher dimensions.
\end{abstract}

\section{Introduction}

A finitely generated group $G$ is said to be \emph{hyperbolic}, if the Cayley graph of $G$ with respect to a finite generating set $S$ has \emph{slim triangles}; that is, there exists some number $\delta \geq 0$ such that for any triangle $\tau$, each side of $\tau$ is contained in the $\delta$-neighbourhood of the remaining two sides. Since their introduction by Gromov in \cite{Gromov1987}, hyperbolic groups have been of fundamental importance in geometric group theory, providing a template to study wider classes of groups \cite{Farb1998,  Osin2016, BehrstockHagenSisto2017}.

It is well known that the class of hyperbolic groups is not closed with respect to taking subgroups. Indeed, the free group of rank two $F_2$ contains a copy of the free group of infinite rank $F_{\infty}$, which is not hyperbolic since it is not finitely generated. Note that every torsion-free hyperbolic group admits a finite classifying space given by the Rips complex \cite[p.468]{BridsonHaefliger1999}. Hence, a natural question to ask is whether the only obstruction for a subgroup of a hyperbolic group $G$ to be hyperbolic arises from its failure to satisfy sufficiently strong \emph{finiteness properties} (see \cref{sec:criterion} for definitions).

The study of this question has a long history, starting with the work of Rips who constructed the first examples of non-hyperbolic subgroups of hyperbolic groups which are finitely generated \cite{Rips1982}. Then Brady produced a non-hyperbolic subgroup of a hyperbolic group which is finitely presented but not of type $\mathrm{FP}_3(\mathbb{Z})$ \cite{Brady1999}. The existence of non-hyperbolic subgroups of hyperbolic groups with higher finiteness properties was observed by Isenrich--Martelli--Py \cite{IsenrichMartelliPy2021} and Fisher \cite{Fisher2021}. Most recently, Isenrich--Py, using methods from complex geometry, produced non-hyperbolic subgroups of hyperbolic groups of type $\mathrm{F}_n$ but not $\mathrm{F}_{n+1},$ for every $n\geq 1$ \cite[Corollary 3]{IsenrichPy2022}. 

A natural question that follows is whether one can obtain non-hyperbolic subgroups of hyperbolic groups of type $\mathrm{FP}(R)$ for some ring $R$.  Italiano--Martelli--Migliorini recently obtained the first such known example, which is moreover of type $\mathrm{F}$ \cite{ItalianoMartelliMigliorini2023}. To do so, they constructed a topological fibring of a finite-volume hyperbolic 5-manifold $M$; the resulting fibre is a finite-volume 4-manifold $N$ such that the outer automorphism group $\mathrm{Out}\,\pi_1(N)$ is infinite. They argue further that one can fill the cusps in $M$ to obtain a negatively curved pseudo-manifold $M'$ which still fibres, and such that the new fibre $N'$ is aspherical and has $\mathrm{Out}\, \pi_1(N')$ infinite. If $\pi_1(N')$ were to be hyperbolic, Paulin's theorem \cite{Paulin1991} would imply that $\pi_1(N')$ splits over an infinite cyclic subgroup. However, a simple Mayer--Vietoris argument shows that $\pi_1(N')$ cannot split over $\Z$.

In the first part of this note, we give a general criterion for constructing non-hyperbolic subgroups of hyperbolic groups which satisfy strong finiteness properties. A group $G$ is said to be \emph{special} if it is the fundamental group of a special cube complex in the sense of Haglund--Wise \cite{HaglundWise2008}.

\begin{thmx}\label{thm_criterion}
    Let $G$ be a torsion-free hyperbolic virtually special group which is $L^2$-acyclic and such that $\mathrm{cd}_{\mathbb{Q}}(G) \geq 4$. Then $G$ contains a non-hyperbolic subgroup $N \leq G$ of type $\mathrm{FP}(\mathbb{Q})$.
\end{thmx}

A group $G$ is \emph{$L^2$-acyclic} if its $L^2$-Betti numbers $b_i^{(2)}(G)$ vanish in every degree. The definition of $L^2$-Betti numbers is technical and not necessary for any of the arguments in this note, and thus it will be omitted. We refer the interested reader to L\"{u}ck's excellent book on the subject \cite{Lueck2002}.

The proof of the criterion crucially uses the main result of Fisher in \cite{Fisher2021} which shows how to obtain homomorphisms to $\Z$ with kernels that have strong finiteness properties. We note that it was already observed by Fisher \cite{Fisher2021} and Isenrich--Martelli--Py \cite{IsenrichMartelliPy2021} that such an approach can be used to construct non-hyperbolic subgroups of hyperbolic groups of type $\mathrm{FP}_n(\mathbb{Q})$ but not $\mathrm{FP}_{n+1}(\mathbb{Q})$. 

The main application of the criterion in \cref{thm_criterion} is the construction of infinitely many non-hyperbolic subgroups of hyperbolic groups which are of type $\mathrm{FP}(\mathbb{Q})$:

\medskip
\begin{duplicate}[\cref{examples}]
    For every integer $n\geq 2$, there exists a non-hyperbolic subgroup $N$ of a torsion-free hyperbolic group, such that $N$ is of type $\mathrm{FP}(\mathbb{Q})$ and $\mathrm{cd}_{\mathbb{Q}}(N) = 2n$. 
\end{duplicate}
\medskip

Let $\mathcal{F}$ denote the family of discrete countable groups with finite cohomological dimension over $\mathbb{Z}$. By the work of Sauer \cite[Theore~1.2]{Sauer2006}, cohomological dimension over $\mathbb{Z}$ is a quasi-isometry invariant amongst groups in $\mathcal{F}$. Hence we get the following corollary: 

\begin{corx}\label{corollary}
    There exists infinitely many quasi-isometry classes of finitely generated subgroups of hyperbolic groups which are of type $\mathrm{FP}(\mathbb{Q})$ and which are not hyperbolic. 
\end{corx}

Recall that a group $G$ is said to be \emph{locally hyperbolic} if every finitely generated subgroup of $G$ is hyperbolic. As another application of our methods, we show that local hyperbolicity and algebraic fibring are mutually exclusive phenomena in groups of higher cohomological dimension: 

\begin{corx}\label{cor_locally_hyperbolic}
    Let $G$ be a finitely generated, locally hyperbolic, torsion-free group. Suppose that $G$ fibres algebraically. Then 
    \[G \simeq (F_n \ast \Sigma_1 \ast \ldots \ast \Sigma_k) \rtimes \mathbb{Z},\]
    where $F_n$ is a free group of rank $n$, and each $\Sigma_i$ is a closed surface group. In particular, $\mathrm{cd}_{\mathbb{Z}}(G) \leq 3$.
    
    Furthermore, if $\mathrm{cd}_{\mathbb{Q}}(G) = 3$ then $G$ contains a closed surface subgroup.
\end{corx}

 \subsection*{Acknowledgements}

 The author would like to thank her supervisors Martin Bridson and Dawid Kielak for their invaluable support and guidance, and Martin for pointing out the examples of arithmetic lattices of simplest type. Thanks also go to Marco Linton for suggesting a more general version of \cref{main}, and the anonymous referee for their comments which greatly improved the exposition in this paper.

 This work has received funding from an Engineering and Physical Sciences Research Council studentship (Project Reference 2422910). 

\section{Preliminaries}\label{sec:preliminaries}

Let $R$ be a ring. A group $G$ is said to be \emph{of type $\mathrm{FP}_n(R)$} if there exists a projective resolution
\[\cdots \to P_{n+1} \to P_n \to \cdots \to P_0 \to M \to 0,\]
where each $P_i$ is a left $RG$-module which is finitely generated whenever $i \leq n,$ and $M$ is the trivial $RG$-module. We say that $G$ is \emph{of type $\mathrm{FP}(R)$} if there exists a finite projective resolution of $M$ by finitely generated $RG$-modules. 

A group $G$ \emph{is of type $F_n$} if it has a classifying space with a finite $n$-skeleton, and of type $F$ if it has a finite classifying space. We note that a group $G$ is of type $\mathrm{FP}_1(R)$ over any ring $R$ if and only if it is of type $F_1$, and the latter condition is equivalent to $G$ being finitely generated. Any group of type $F_n$ is also of type $\mathrm{FP}_n(R)$ for any ring $R$, however the converse does not hold in general (see \cite[Section~6]{BestvinaBrady1997}).

\begin{defn}
   The \emph{cohomological dimension} $\mathrm{cd}_R(G)$ of a group $G$ over a ring $R$ is the minimal integer $d$ such that the trivial $RG$-module $R$ has a projective resolution of length $d$ by $RG$-modules, provided that such an integer $d$ exists. Otherwise we set $\mathrm{cd}_R(G) = \infty$.
\end{defn}

Consider a short exact sequence of groups 
    \[1 \to N \to G \to Q \to 1.\]
The Lyndon--Hoschild--Serre spectral sequence for a short exact sequence of groups implies that for any ring $R$, the cohomological dimensions of the groups in the short exact sequence over $R$ satisfy
    \begin{equation}\label{eq:sesincd}\mathrm{cd}_{R}(G) \leq \mathrm{cd}_{R}(N) + \mathrm{cd}_{R}(Q).\end{equation}
There are specific circumstance in which the inequality in \eqref{eq:sesincd} can be turned into an equality:

\begin{thm}[Fel'dman \cite{Feldman1971}]\label{Feldman}
    Let $R$ be a ring. Consider a short exact sequence of groups 
    \[1 \to N \to G \to Q \to 1.\] Suppose that $N$ is of type $\mathrm{FP}(R)$ and $H^n(N ; RN)$ is $R$-free, where $n = \mathrm{cd}_R(N)$. 
    
    If $\mathrm{cd}_R(G) < \infty$ then 
    \[\mathrm{cd}_{R}(G) = \mathrm{cd}_{R}(N) + \mathrm{cd}_{R}(Q).\]
\end{thm}
\smallskip

\begin{remark}
    If $R = \mathbb{Q}$ then $H^n(N; RN)$ is $R$-free.
\end{remark}

Let $X$ be a locally finite graph. For any finite collection $K$ of edges of $X$, let $e(X; K)$ be the number of infinite connected components of the graph $X \setminus K$. The \emph{number of ends of $X$} is defined to be \[e(X) :=\mathrm{sup}\ e(X; K),\]
where the supremum is taken over all finite collections of edges of $X$.

\begin{defn}
 If $G$ is a finitely generated group then the \emph{number of ends of $G$} is the number of ends of a Cayley graph $\mathrm{Cay}(G,S)$ of $G$ with respect to any finite generating set $S$ of $G$. 
\end{defn} 

The number of ends of a finitely generated group $G$ satisfies $e(G) \in \{0, 1,2, \infty \}$ and $e(G) = 2$ exactly when $G$ is virtually infinite cyclic \cite{Hopf1944, Freudenthal1945}. Moreover, by the work of Stallings \cite{Stallings1968}, if $G$ is finitely generated and torsion free, then $e(G) = \infty$ if and only if $G$ splits as a free product $G = G_1 \ast G_2$ with $G_1$ and $G_2$ non-trivial.

\section{A criterion for exotic subgroups of hyperbolic groups}\label{sec:criterion}

The aim of this section is to prove a general criterion for constructing exotic subgroups of hyperbolic groups: 

%We write $\mathrm{cd}(G)$ to denote the cohomological dimension of $G$ over $\mathbb{Z}$. 
\smallskip

\begin{thm}\label{thm_criterion_improved}
    Let $G$ be a torsion-free hyperbolic virtually special group with $\mathrm{cd}_{\mathbb{Q}}(G) \geq 4$. Suppose that the $L^2$-Betti numbers of $G$ satisfy $b_{i}^{(2)}(G) = 0$ for all $i \leq n.$ Then $G$ contains a non-hyperbolic subgroup $N \leq G$ of type $\mathrm{FP}_n(\mathbb{Q})$.
    
    Moreover, $\mathrm{cd}_{\mathbb{Q}}(N) \in \{\mathrm{cd}_{\mathbb{Q}}(G) -1, \mathrm{cd}_{\mathbb{Q}}(G) \}$ and if $\mathrm{cd}_{\mathbb{Q}}(G) \leq n$ then $\mathrm{cd}_{\mathbb{Q}}(N) =  \mathrm{cd}_{\mathbb{Q}}(G) -1$.
\end{thm}
\smallskip

Note that if a group $N$ is of type $\mathrm{FP}_n(\mathbb{Q})$ and its cohomological dimension satisfies $\mathrm{cd}_{\mathbb{Q}}(N) \leq n$, then $N$ is of type $\mathrm{FP}(\mathbb{Q})$. Hence \cref{thm_criterion_improved} implies \cref{thm_criterion} from the introduction.

We delay the proof of \cref{thm_criterion_improved} to \cref{sec:proof}. First in \cref{sec:hypKernels} we study the structure of hyperbolic kernels of infinite quotients $G \twoheadrightarrow Q$, where $G$ is a hyperbolic group. As an application, we prove \cref{cor_locally_hyperbolic} on the structure of locally hyperbolic groups which fibre.

\subsection{Hyperbolic subgroups of hyperbolic groups}\label{sec:hypKernels}

\begin{lemma}\label{one-ended} Let $H$ be a one-ended torsion-free hyperbolic group. If $H$ is not the fundamental group of a surface then $H \rtimes \mathbb{Z}$ contains a subgroup isomorphic to $\mathbb{Z}^2$. 
\end{lemma}

\begin{proof}

    Suppose that $H$ is one-ended torsion-free hyperbolic and that it is not the fundamental group of a closed surface. Paulin's theorem \cite{Paulin1991} (see also \cite{BridsonSwarup1994}) shows that if $\mathrm{Out}(H)$ is infinite then $H$ acts on an $\mathbb{R}$-tree with cyclic arc stabilisers. It then follows from the work of Rips that $H$ admits such an action on a simplicial tree, i.e. it splits over $\Z$; see \cite{BestvinaFeighn1995}. 
    
    Hence if $H$ does not split over $\Z$ then $\mathrm{Out}(H)$ is finite. In particular, there exists $m \geq 1$ such that $\phi^m$ is an inner automorphism of $H$. It follows that there exists some $x\in G$ such that $\langle H, xt^m\rangle \simeq H  \times \mathbb{Z}$, and so $G$ contains a subgroup isomorphic to $\Z^2$.
    
   Suppose now that $H$ splits over $\mathbb{Z}$. Since $H$ is not a surface group, it admits a non-trivial JSJ decomposition \cite[Theorem~1.7]{Sela1997}. The JSJ splitting is canonical, and thus $\phi$ permutes the conjugacy classes of edge stabilisers. Hence, for every edge of the splitting, there exists some $x \in G$ and $m \geq 1$ such that the generator of the edge stabiliser commutes with $t^mx$, and thus $G$ contains a $\Z^2$ subgroup.
\end{proof}

\begin{prop}\label{main} Let $G$ be a group with no $\mathbb{Z}^2$-subgroups. Suppose that $G$ fits into the short exact sequence of groups
\[1 \to H \to G \to Q \to 1,\]
where $Q$ is infinite. If $H$ is torsion-free hyperbolic, then $H$ splits as a free product of finitely many fundamental groups of compact hyperbolic surfaces.

In particular, the cohomological dimension of $H$ over $\Z$ satisfies $\mathrm{cd}_{\mathbb{Z}}(H)  \leq 2$.
\end{prop}

\begin{proof}
The quotient $Q$ is hyperbolic and thus contains an element of infinite order by \cite[Theorem~A]{Mosher1996}. Let $\gamma$ denote the lift of such an element to $G$. Then $G$ contains a subgroup $\langle H , \gamma \rangle \simeq H \rtimes \Z$. If $H$ is one-ended, then by Lemma~\ref{one-ended} $H$ is the fundamental group of a closed hyperbolic surface and $\mathrm{cd}_{\Z}(H) = 2$.
 
 If $H$ is two-ended then it is isomorphic to the infinite cyclic group and $G \simeq \Z \rtimes \Z$. But then $G$ is virtually $\Z^2$ and thus it cannot be hyperbolic. 

Suppose that $H$ is infinitely-ended. Let $H = H_1 \ast \ldots \ast H_n \ast F_k$ be the Grushko decomposition of $H$. Fix a lift $t\in G$ of a generator of the infinite cyclic group. Let $\phi \in \mathrm{Aut}(H)$ be the automorphism of $H$ induced by the conjugation action of $t$ on $H$ in $G$. Then $\phi$ permutes the conjugacy classes of the subgroups $H_i$. Hence for all $i$, there exists $x_i \in G$ such that \[\langle H_i, t^mx_i \rangle \simeq H_i \rtimes_{\mathrm{Ad}_{x_i} \phi^m} \mathbb{Z},\] 
where $\mathrm{Ad}_{x_i}$ denotes the inner automorphism of $G$ which acts by conjugation with $x_i$. Note that each $H_i$ is one-ended, torsion-free hyperbolic. If there exists some $H_i$ which is not a surface group, then by Lemma~\ref{one-ended} $G$ is not hyperbolic. Moreover, if $H = H_1 \ast \ldots \ast H_n \ast F_k$ then the cohomological dimension of $H$ satisfies
\[\mathrm{cd}_{\Z}(G) \leq \mathrm{max}\{ \mathrm{cd}_{\mathbb{Z}}(H_i), \mathrm{cd}_{\mathbb{Z}}(F_k)\} \leq 2.\]\end{proof}

We note that the special case of \cref{main} when $Q \simeq \Z$ was observed by Peter Brinkmann in his thesis \cite{BrinkmannThesis2000}. The author would like to thank Daniel Groves for pointing this out. 
\smallskip

We end this subsection by outlining the proof of \cref{cor_locally_hyperbolic}. Let $G$ be a locally hyperbolic torsion-free group which fibres algebraically. Then the kernel $H$ of the fibration of $G$ is hyperbolic and thus by \cref{main}, it is a free product of closed surface groups and a free group. Arguing as in the proof of \cref{main}, $\mathrm{cd}_{\Z}(H) \leq 2$ and by \eqref{eq:sesincd},
\[\mathrm{cd}_{\Z}(G) \leq \mathrm{cd}_{\Z}(H) + \mathrm{cd}_{\Z}(\Z) \leq 3.\] This proves the first part of \cref{cor_locally_hyperbolic}. 

For the second part, note that $H$ is of type $\mathrm{FP}(\mathbb{Q})$ and thus by Fel'dman's \cref{Feldman}, $\mathrm{cd}_{\mathbb{Q}}(G) = \mathrm{cd}_{\mathbb{Q}}(H) + 1$.  The cohomological dimension of $H$ over $\mathbb{Q}$ is bounded above by 2, and it is exactly equal to 2 when at least one of the components of the free product decomposition of $H$ is a closed surface group. This proves the second part of \cref{cor_locally_hyperbolic}.

\subsection{Proof of the criterion}\label{sec:proof}

The key tool for proving the criterion is the following theorem of Fisher:

\begin{thm}[Fisher \cite{Fisher2021}]\label{thm_Fisher}
    Let $G$ be a group of type $\mathrm{FP}_n(\mathbb{Q})$. Suppose that $G$ is RFRS and its $L^2$-Betti numbers satisfy $b_{i}^{(2)}(G) = 0$ for $i \leq n$. Then there exists a finite index subroup $G' \leq G$ and an epimorphism $\varphi \colon G' \to \Z$ such that $\ker \varphi$ is of type $\mathrm{FP}_n(\mathbb{Q})$.
\end{thm}

\begin{proof}[Proof of \cref{thm_criterion_improved}]
    If $G$ is torsion-free hyperbolic, then $G$ has a finite classifying space given by the Rips complex (see \cite[p.468]{BridsonHaefliger1999}). In particular, $G$ is of type $\mathrm{FP}_{n}(\mathbb{Q})$ for all $n \geq 1$. Hence by Fisher's  \cref{thm_Fisher}, the hypothesis on the vanishing of the $L^2$-Betti numbers of $G$ implies that there exists a finite index subgroup $G' \leq G$ such that $G'$ fibres algebraically with kernel $N$ of type $\mathrm{FP}_k(\mathbb{Q})$ for $k \leq n$. 
    
    Consider the Lyndon--Hoschild--Serre spectral sequence for the short exact sequence  
    \[1 \to N \to G' \to G'/N \to 1. \]  It follows that
    \[\mathrm{cd}_{\mathbb{Q}}(G') \leq \mathrm{cd}_{\mathbb{Q}}(N) + \mathrm{cd}_{\mathbb{Q}}(G'/N).\]
    \smallskip
Since $G' / N \simeq \Z$, we have that $\mathrm{cd}_{\mathbb{Q}}(G' / N) = 1$. Moreover since $G' \leq G$ is a finite index subgroup, \cite[Proposition~5.7]{Bieri1981} implies that the cohomological dimensions of $G$ and $G'$ over $\mathbb{Q}$ agree.

Combining these facts, we obtain
  \begin{equation}\label{eq:1}\begin{split} \mathrm{cd}_{\mathbb{Q}}(N) &\geq \mathrm{cd}_{\mathbb{Q}}(G') - 1  \\
    &= \mathrm{cd}_{\mathbb{Q}}(G) - 1 \geq 3. \end{split}\end{equation} 
Hence $\mathrm{cd}_{\Z}(N) \geq \mathrm{cd}_{\mathbb{Q}}(N) \geq 3$ and thus by \cref{main} the subgroup $N$ is not hyperbolic. 

If $\mathrm{cd}_{\mathbb{Q}}(G)$ is bounded above by $n$, then 
\[\mathrm{cd}_{\mathbb{Q}}(N)  \leq \mathrm{cd}_{\mathbb{Q}}(G)\leq n.\]
Since $N$ is of type $FP_n(\mathbb{Q})$ and satisfies $\mathrm{cd}_{\mathbb{Q}}(N) \leq n$, we conclude that $N$ is of type $FP(\mathbb{Q})$. Thus by Fel'dman's \cref{Feldman} the first inequality in \eqref{eq:1} is an equality. \end{proof}

\section{Explicit examples}

In this section we will construct explicit examples of non-hyperbolic subgroups of hyperbolic groups with strong finiteness properties.

\smallskip
\begin{thmx}\label{examples}
    For every integer $n\geq 2$, there exists a non-hyperbolic subgroup $N$ of a torsion-free hyperbolic group, such that $N$ is of type $\mathrm{FP}(\mathbb{Q})$ and $\mathrm{cd}_{\mathbb{Q}}(N) = 2n$. 
\end{thmx}
\smallskip

The non-hyperbolic subgroups in \cref{examples} will arise as subgroups of fundamental groups of hyperbolic manifolds. Dodziuk shows that it is possible to extract information about the $L^2$-cohomology of a manifold $M$ from studying the $L^2$-harmonic forms on $M$ \cite{Dodziuk1977}. His explicit calculation of the $L^2$-harmonic forms for hyperbolic manifolds yields the following result about the $L^2$-Betti numbers of such manifolds.

\begin{thm}[Dodziuk \cite{Dodziuk1979}] \label{Dodziuk}
Let $n \geq 1$ be odd. If $G$ is the fundamental group of a closed oriented hyperbolic $n$-manifold, then for all $i \geq 0$, 
\[b_i^{(2)}(G) = 0.\]
\end{thm}

Our criterion in \cref{thm_criterion_improved} combined with Dodziuk's \cref{Dodziuk} shows that any odd-dimensional closed hyperbolic $n$-manifold with $n \geq 4$ whose fundamental group is virtually special, will contain a subgroup as in Theorem~\ref{examples}. 

Such examples can be obtained from the work of Bergeron--Haglund--Wise, who prove that in every dimension $n$ there exist cocompact arithmetic lattices in $SO(n,1)$ which are virtually special \cite[Theorem~1.10]{BergeronHaglundWise2011}. After passing to a subgroup of finite index, we may assume that the lattice $G$ is torsion-free and special. 

Note that in this case there is an alternate way to see that the kernel $K$ of an algebraic fibration is non-hyperbolic. Mainly, if $K$ is of type $FP_{n-1}(\mathbb{Q})$ then by a theorem of Hillman \cite[Theorem~1.19]{Hillman2002} it is an $(n-1)$-dimensional Poincar\'{e} Duality group  (over $\mathbb{Q}$). Thus, a Mayer--Vietoris argument shows that it cannot split over $\Z$ or over the trivial group. Assuming that $K$ is hyperbolic, it follows by Paulin's theorem \cite{Paulin1991} and the Rips machine \cite{BestvinaFeighn1995} that the outer automorphism group of $K$ is finite. Then, arguing as in the proof of Lemma~\ref{one-ended}, we conclude that $G$ virtually splits as a direct product $H \times \Z$ and thus cannot be hyperbolic.

Finally, we note that it would be interesting to construct such examples which do not arise from the world of manifolds. One potential avenue to do so is through the work of Arenas in \cite{Arenas2022}, whose cubical Rips construction is a new tool for constructing higher dimensional Gromov hyperbolic groups which fibre algebraically.

\bibliographystyle{alpha}
\bibliography{main.bib}

\end{document}